\documentclass[a4paper,12pt,dvipdfmx,pdftex]{article}

\pdfoutput=1
\title{
The Riemann hypothesis and holomorphic index \\
in complex dynamics
}
\author{
Tomoki Kawahira
\thanks{Research partially supported by JSPS}\\
Tokyo Institute of Technology
}

\usepackage{enumerate}
\usepackage[dvipdfmx]{color}
\usepackage[dvipdfmx]{graphicx}
\usepackage{amsmath}
\usepackage{amscd}
\usepackage{amssymb}
\usepackage{amsfonts}
\usepackage{theorem}

\newtheorem{thm}{Theorem}
\newtheorem{prop}[thm]{Proposition}

\newtheorem{cor}[thm]{Corollary}
\theorembodyfont{\rmfamily}

\newcommand{\thmref}[1]{Theorem \ref{#1}}
\newcommand{\propref}[1]{Proposition \ref{#1}}
\newcommand{\corref}[1]{Corollary \ref{#1}}

\newcommand{\figref}[1]{Figure \ref{#1}}

\renewcommand{\Bar}[1]{\overline{#1}}

\newcommand{\parag}[1]{
\medskip
\noindent {\bf #1}
}

\newcommand{\C}{\ensuremath{\mathbb{C}}}

\newcommand{\Q}{\ensuremath{\mathbb{Q}}}

\newcommand{\Chat}{\ensuremath{\widehat{\mathbb{C}}}}
\newcommand{\R}{\ensuremath{\mathbb{R}}}

\newcommand{\D}{\ensuremath{\mathbb{D}}}

\newcommand{\Z}{\ensuremath{\mathbb{Z}}}
\newcommand{\N}{\ensuremath{\mathbb{N}}}

\newcommand{\kakko}[1]{{\left( #1 \right)}}

\newcommand{\skakko}[1]{{\left\{ #1 \right\}}}

\renewcommand{\Re}{\mathrm{Re}\,}
\renewcommand{\Im}{\mathrm{Im}\,}

\newcommand{\id}{\mathrm{id}}

\newcommand{\al}{{\alpha}}
\newcommand{\gam}{{\gamma}}
\newcommand{\lam}{{\lambda}}

\newcommand{\zbar}{\Bar{z}}

\newcommand{\cc}{\circ}
\newcommand{\dz}{\, dz}

\newcommand{\mapsTo}[1]{\stackrel{{#1}}{\longmapsto}}

\newcommand{\ee}{~=~}

\newcommand{\st}{\,:\,}

\newcommand{\QED}{\hfill $\blacksquare$}
\renewcommand{\Bar}{\overline}

\begin{document}
\maketitle

\begin{abstract}
We give an interpretation of the Riemann hypothesis 
in terms of complex and topological dynamics.
For example, 
the Riemann hypothesis is affirmative and 
all zeros of the Riemann zeta function are simple
if and only if 
a certain meromorphic function has 
no attracting fixed point.
To obtain this, we use holomorphic index 
(residue fixed point index), 
which characterizes local properties of 
fixed points in complex dynamics.
\end{abstract}


\section{The Riemann zeta function}
For $s \in \C$, the series
$$
\zeta(s) \ee 1+\frac{1}{2^s}+\frac{1}{3^s} +\cdots 
$$
converges if $\Re s >1$. 
Indeed, $\zeta(s)$ is analytic on the half-plane
$\skakko{s \in \C \st \Re s  > 1}$ 
and continued analytically to 
a meromorphic function 
(a holomorphic map) $\zeta: \C \to \Chat = \C \cup \skakko{\infty}$ 
with only one pole at $s=1$, which is simple. 
This is \textit{the Riemann zeta function}.

It is known that $\zeta(s) = 0$ when $s = -2,-4,-6,\ldots$. 
These zeros are called {\it trivial zeros} of 
the Riemann zeta function. We say the other zeros are {\it non-trivial}.

The {\it Riemann hypothesis}, 
which is the most important conjecture on the Riemann zeta function,
concerns the alignment of the non-trivial zeros of $\zeta$: 

\parag{The Riemann hypothesis.}
{\it All non-trivial zeros lie on the vertical line 
$\{s \in \C \st \Re s = 1/2\}$.}\\

The line $\{s \in \C \st \Re s = 1/2\}$ is called {\it the critical line}. 
It is numerically verified that the first 
$10^{13}$ zeros (from below) lie on the critical line \cite{G}. 

It is also conjectured that every zero of $\zeta$ is simple. 
We refer to this conjecture as {\it the simplicity hypothesis} 
after some literature. (See \cite{RS} or \cite{Mu} for example. 
One may find some related results and observations 
in \cite[\S 10.29, \S 14.34, \S 14.36]{T}. )

The aim of this note is to translate the Riemann hypothesis 
in terms of complex and topological dynamical systems.
More precisely, we translate the locations of the non-trivial zeros 
into some dynamical properties of the fixed points of a meromorphic function of the form
$$
\nu_g(z) = z-\frac{g(z)}{z\,g'(z)},
$$
where $g$ is a meromorphic function on $\C$ that shares (non-trivial) zeros with $\zeta$. 
For example, we will set $g = \zeta$ or the {\it Riemann xi function}
$$
\xi(z) : =  \frac{1}{2}z(1-z)\pi^{-z/2}
\Gamma\kakko{\frac{z}{2}}\,\zeta(z),
$$
etc. The function $\nu_g(z)$ is carefully chosen so that 
\begin{itemize}
\item
If $g(\al) = 0$ then $\nu_g(\al) = \al$; and
\item
The {\it holomorphic index} (or the {\it residue fixed point index}, 
see \S 2) of $\nu_g$ at $\al$ is $\al$ itself 
when $\al$ is a simple zero of $g(z)$.
\end{itemize}
See \S 3 for more details.

For a given meromorphic function $g:\C \to \Chat$,
we say a fixed point $\al$ of $g$ is
{\it attracting} if $|g'(\al)|<1$, 
{\it indifferent} if $|g'(\al)| = 1$, and 
{\it repelling} if $|g'(\al)|>1$.
Here are some possible translations of the Riemann hypothesis
 (plus the simplicity hypothesis) by $\nu_\zeta$:
 
\begin{thm}\label{thm_zeta}
The following conditions are equivalent:
\begin{enumerate}[\rm (a)]
\item
The Riemann hypothesis is affirmative and every non-trivial zero of $\zeta$ is simple.
\item 
Every non-trivial zero of $\zeta$ is an indifferent fixed point of 
the meromorphic function $\nu_\zeta(z) := z-\dfrac{\zeta(z)}{z\,\zeta'(z)}$ .
\item 
The meromorphic function $\nu_\zeta$ above has no attracting fixed point.
\item 
There is no topological disk $D$ with $\Bar{\nu_\zeta(D)} \subset D$.
\end{enumerate}
\end{thm}

Note that (d) is a topological property of the map $\nu_\zeta:\C \to \Chat$,
in contrast to the analytic (or geometric) nature of (a).
We will present some variants of this theorem in \S 4. 

\S 5 is devoted for some numerical observations
 and questions on linearization problem. 
\S 6 is an appendix: we apply Newton's method to the Riemann zeta function. 
We will also give a ``semi-topological" criteria for the Riemann hypothesis 
in terms of the Newton map (Theorem \ref{thm_NZeta}).

\parag{Remarks.}

\begin{itemize}
\item
We can apply the method of this note to the $L$-functions without extra effort.
\item
The following well-known facts are implicitly used in this paper:
\begin{itemize}
\item
Every non-trivial zero is located in the {\it critical stripe} 
$$
\mathcal{S}:=\skakko{s \in \C \st 0 < \Re s < 1}.
$$
\item
The functional equation
$\zeta(s) = 2^s \pi^{s-1}\sin(\pi s/2)\Gamma(1-s)\zeta(1-s)$
implies that if $\al$ is a non-trivial zero of $\zeta(s)$, 
then so is $1-\al$ and they have the same order.
\end{itemize} 
Hence the non-trivial zeros are symmetrically arrayed with respect to $s = 1/2$.
By these properties, we will mainly consider the zeros 
which lie on the upper half of the critical stripe $\mathcal{S}$.
See \cite{E, T} for more details.
\item
We used {\it Mathematica} 10.0 for all the numerical calculation. 
\end{itemize}

\parag{Acknowledgments.}
The author would like to thank Masatoshi Suzuki for helpful comments.

\section{Fixed points and holomorphic indices}

\parag{Multiplier.}
Let $g$ be a holomorphic function on a domain $\Omega \subset \C$.
We say $\al \in \Omega$ is a {\it fixed point} of $g$ with 
{\it multiplier} $\lam \in \C$ if $g(\al)=\al$ and $g'(\al)=\lam$.
The multiplier $\lam$ is the primary factor that determines the local dynamics near $\al$. 
In fact, the Taylor expansion about $\al$ gives a 
representation of the local action of $g$:
\begin{equation}\label{eq_Taylor}
g(z) -\al \ee \lam(z-\al) +O(|z-\al|^2).
\end{equation}
Hence the action of $g$ near $\al$
is locally approximated by $w \mapsto \lam w$ in the coordinate $w = z-\al$.
We say the fixed point $\al$ is 
{\it attracting, indifferent}, or {\it repelling}
according to $|\lam|<1$, $|\lam|=1$, or $|\lam|>1$.

\parag{Topological characterization.}
Attracting and repelling fixed points of holomorphic mappings 
have purely topological characterizations (cf. Milnor \cite[\S 8]{Mi}):

\begin{prop}[Topological characterization of fixed points]
\label{prop_att_rep}
Let $g$ be a holomorphic function on a domain $\Omega \subset \C$.
The function $g$ has an attracting (resp. repelling) fixed point  
if and only if 
there exists a topological disk $D \subset \Omega$
such that $\overline{g(D)} \subset D$
(resp. $g|_D$ is injective and $\Bar{D}\subset g(D) \subset \Omega$).
\end{prop}
The condition that $g$ is holomorphic is essential. 
For example, the proposition is false if 
we only assume that $g$ is $C^\infty$.

\begin{figure}[htbp]
\begin{center}
\includegraphics[width=0.5\textwidth]{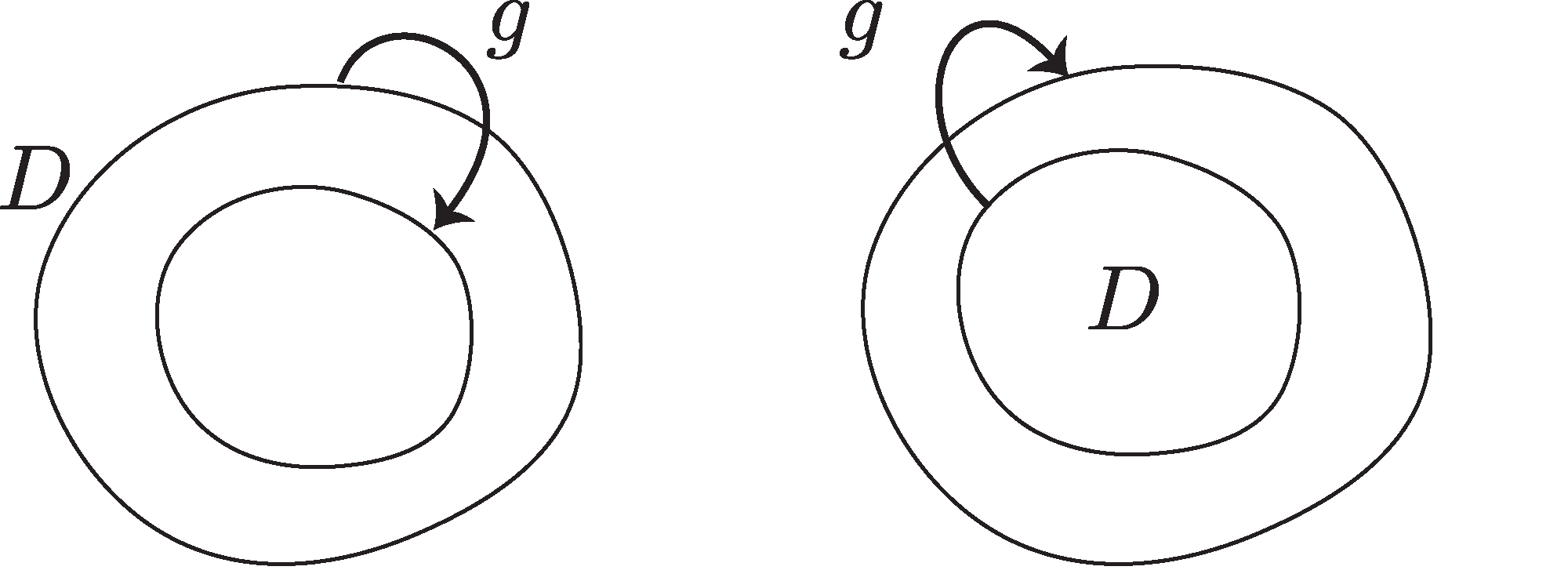}
\end{center}
\caption{Topological disks which contain an attracting or a repelling fixed point.}
\label{fig_}
\end{figure}

\parag{Proof.}
Suppose that $g$ has an attracting fixed point $\al \in \Omega$.
By the Taylor expansion (\ref{eq_Taylor}) as above,
we have $g(z)-\al=(\lam+o(1))(z-\al)$ near $\al$ 
and thus a small circle around $\al$ is mapped 
strictly inside the circle. 
By the maximum principle, this circle bounds a round disk $D$ such that 
$\overline{g(D)} \subset D$.

Conversely, if a topological disk $D$ in $\Omega$ 
satisfies $\overline{g(D)} \subset D$, 
we may observe the map $g:D \to g(D) \subset D$ via the Riemann
map and it is enough to consider the case of $D=\D$, the unit disk.
By the Schwarz-Pick theorem (see \cite[\S 1]{Ah2}) the map is 
strictly contracting with respect to the distance 
$d(z,w):=|z-w|/|1-\zbar w|$ on $\D$. 
Hence there exists a fixed point $\al$ by the fixed point theorem. 
It must be an attracting fixed point by the Schwarz lemma.

The repelling case is analogous.
\QED

\parag{Holomorphic index.}
Let $\al$ be a fixed point of 
a holomorphic function $g:\Omega \to \C$.
We define the {\it holomorphic index} (or {\it residue fixed point index})
of $\al$ by
$$
\iota(g,\al) := \frac{1}{2 \pi i}\int_C \frac{1}{z-g(z)} \dz,
$$
where $C$ is a small circle around $\al$ with counterclockwise direction. 
The holomorphic index is mostly determined by the multiplier:

\begin{prop}\label{prop_main_formula}
If the multiplier $\lam:=g'(\al)$ is not $1$, then we have $\iota(g, \al) = \dfrac{1}{1-\lambda}.$
\end{prop}
See \cite[Lem.12.2]{Mi} for the proof.

\parag{Remark.}
Any complex number $K$ may be the holomorphic index 
of a fixed point of multiplier $1$.
For example, the polynomial $g(z)=z-z^2+Kz^3$ has a 
fixed point at zero with $g'(0)=1$ and $\iota(g,0)=K$.

Since the M\"obius transformation $\lam \mapsto \dfrac{1}{1-\lam}=\iota$ 
sends the unit disk to the half-plane $\skakko{\iota \in \C \st \Re \iota > 1/2}$,
fixed points are classified as follows:

\begin{prop}[Classification by index]
\label{prop_classification}
Suppose that the multiplier $\lam=g'(\al)$ is not $1$. 
Then $\al$ is attracting, repelling, or indifferent
if and only if the holomorphic index $\iota=\iota(g,\al)$ satisfies $\Re \iota >1/2$, 
$<1/2$, or $=1/2$ respectively.
\end{prop}
\begin{center}
\begin{figure}[htbp]
\begin{center}
\includegraphics[width=.7\textwidth]{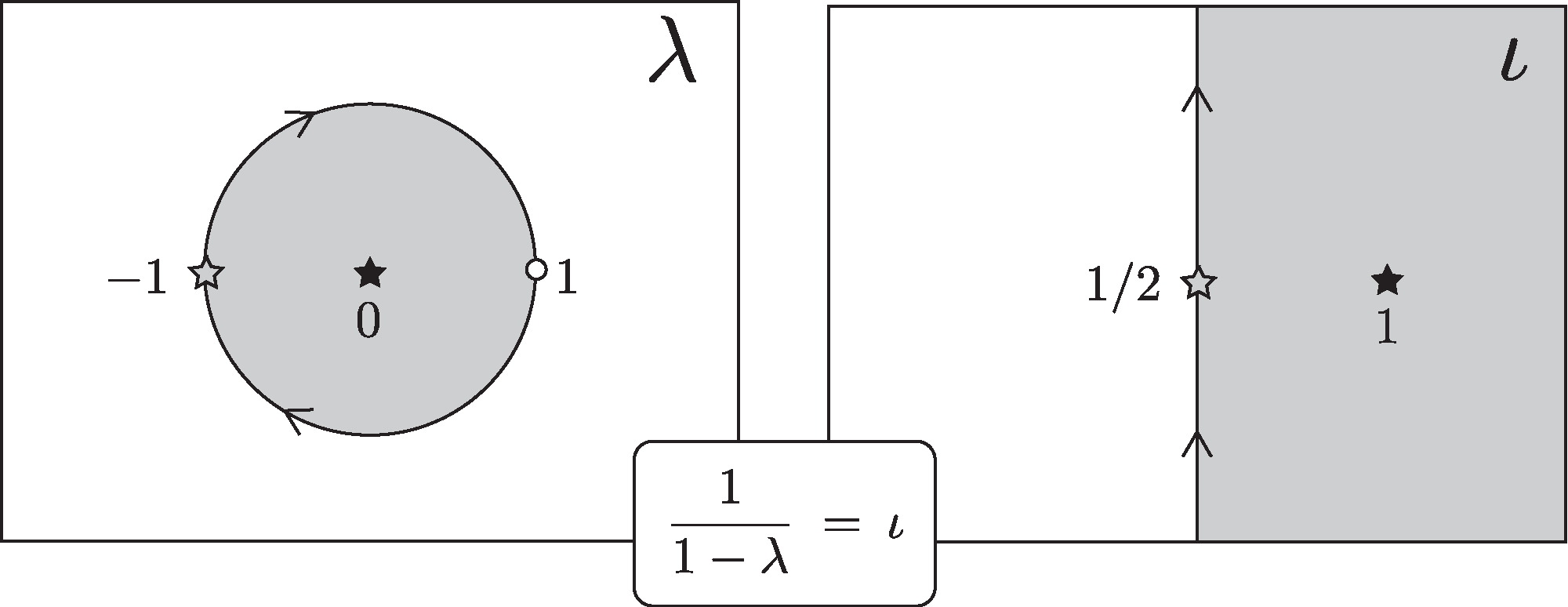}
\end{center}
\caption{Multipliers and holomorphic indeces}
\label{fig_}
\end{figure}
\end{center}
Note that the indifferent fixed points whose multiplier is not $1$ 
correspond to the ``critical line" in the $\iota$-plane.

\section{The nu function}
Let $g:\C \to \Chat$ be a non-constant meromorphic function. 
(We regard such a meromorphic function as 
a holomorphic map onto the Riemann sphere.) 
We define the {\it nu function} $\nu_g:\C \to \Chat$ of $g$ by
$$
\nu_g(z):=z-\frac{g(z)}{z g'(z)}.
$$
This is also a non-constant meromorphic function on $\C$.
We first check:

\begin{prop}[Fixed points of $\nu_g$]
\label{prop_nu_index}
Suppose that $\al \neq 0$. 
Then 
$\al$ is a fixed point of $\nu_g$
if and only if 
$\al$ is a zero or a pole of $g$. 
Moreover, 
\begin{itemize}
\item
if $\al$ is a zero of $g$ 
of order $m \ge 1$,
then $\nu_g'(\al) =1 - \dfrac{1}{m\al}$ 
and $\iota(\nu_g,\al)= m \al$; and
\item 
if $\al$ is a pole of $g$ 
of order $m \ge 1$,
then $\nu_g'(\al) =1 + \dfrac{1}{m\al}$ 
and $\iota(\nu_g,\al)= -m \al$.
\end{itemize}
If $0$ is a zero or a pole of $g$ of order of $m \ge 1$, 
then $\nu_g(0)=1/m$ or $-1/m$ respectively. 
In particular, $0$ is not a fixed point of $\nu_g$. 
\end{prop}
An immediate corollary is:

\begin{cor}
\label{cor_no_fixed_points_with_mult_one}
The function $\nu_g$ has no fixed point 
of multiplier $1$.
\end{cor}
Hence we can always apply 
Proposition \ref{prop_main_formula} to the fixed points of $\nu_g$.

\parag{Proof of Proposition \ref{prop_nu_index}.}
The point $\al \in \C$ is a fixed point 
of $\nu_g(z)=z-g(z)/(zg'(z))$
if and only if $g(\al)/(\al g'(\al))=0$. 
This implies that $\al$ must be a zero or a pole of $g$. 

Conversely, suppose that 
$$
g(z) = C(z-\al)^{M}\kakko{1+O(z-\al)}
$$
for some $M \in \Z-\{0\}$ and $C \in \C-\{0\}$ near $\al$.

Suppose in addition that $\al \neq 0$. 
Then the expansion of $\nu_g$ about $\al$ is
$$
\nu_g(z) = \al + \kakko{1-\frac{1}{M\al}}(z-\al) + O((z-\al)^2)
$$
and thus $\al$ is a fixed point with multiplier
$\nu_g'(\al)=1 - 1/(M\al) \neq 1$.
By Proposition \ref{prop_main_formula},
its holomorphic index is $1/(1-\nu_g'(\al))=M\al$.

If $\al=0$, the expansion of $\nu_g$ about $\al=0$ is
$$
\nu_g(z) =\frac{1}{M} + O(z).
$$
This implies that $\nu_g(0)=1/M \neq 0$.
\QED

\parag{Remark.}
There is another way to calculate the holomorphic index.
(This is actually how the author found the function $\nu_g$.)
Suppose that $\al$ is a zero of $g$ of order $m$.
A variant of the argument principle (\cite[p.153]{Ah1}) 
yields that for any holomorphic function $\phi$ defined near $\al$,
we have 
$$
\frac{1}{2\pi i}\int_C \phi(z) \, \frac{g'(z)}{g(z)} \dz
= m \phi(\al)
$$ 
where $C$ is a small circle around $\al$. 
Set $\phi(z):=z$. Then the equality above is equivalent to 
$$
\frac{1}{2\pi i}\int_C \frac{1}{z-\nu_g(z)}\dz
= m \al. 
$$ 
The same argument is also valid when $\al$ is a pole.

\parag{Example.}
Consider a rational function $g(z)=\dfrac{(z+1)(z-1/2)^2}{z^3(z-1)}$.
By Proposition \ref{prop_classification}
and Proposition \ref{prop_nu_index}, 
the zeros $-1$ and $1/2$ of $g$ are repelling and attracting fixed points 
of $\nu_g$ respectively.
The pole $1$ is a repelling fixed point,
though $0$ is not a fixed point of $\nu_g$.

\section{The Riemann hypothesis}
Let us consider the Riemann zeta function and prove Theorem \ref{thm_zeta}.

\parag{The Riemann zeta.}
Set $g=\zeta$, the Riemann zeta function.
It is known that the trivial zeros 
$\al=-2,\,-4,\,\cdots$ of $\zeta$ are all simple.
By Proposition \ref{prop_nu_index}, 
they are fixed points of $\nu_\zeta$ of multiplier $1-1/\al \neq 1$
and their holomorphic indices are $\al$ itself.
Hence by Proposition \ref{prop_classification}, 
they are all repelling fixed points of $\nu_\zeta$.

Similarly the unique pole $z = 1$ of $\zeta$ is simple and 
it is a repelling fixed point of $\nu_\zeta$ with multiplier $1 + 1/1 = 2$. Hence we have:

\begin{prop}
\label{prop_rep_off_crit_stripe}
Every fixed point $\al$ of $\nu_\zeta$ off the critical stripe 
$\mathcal{S}$ is repelling. 
\end{prop}

Let $\al$ be a non-trivial zero of order $m \ge 1$ 
in the critical strip $\mathcal{S}$. 
(Under the simplicity hypothesis $m$ is always $1$.)
By Proposition \ref{prop_nu_index}, 
$\al$ is a fixed point of $\nu_\zeta$ with multiplier $\lam:=1-1/(m\al) \neq 1$,
and its holomorphic index is $\iota:=m\al$.

If the Riemann hypothesis holds, $\Re \iota=\Re m\al = m/2$.
Thus $\Re \iota=1/2$ if $m=1$ and $\Re \iota \ge 1$ if $m \ge 2$.
Since ``$\Re \iota \ge 1$ in the $\iota$-plane" 
is equivalent to
``$|\lam -1/2|\le 1/2$ (and $\lam \neq 1$) 
in the $\lam$-plane", we have:

\begin{prop}
\label{prop_RH_implies_indiff_or_att}
Under the Riemann hypothesis, 
any fixed point $\al$ of $\nu_\zeta$ in $\mathcal{S}$   
is a zero of $\zeta$ of some order $m \ge 1$ 
that lies on the critical line. Moreover,
\begin{itemize}
\item
when $m=1$, $\al$ is indifferent with multiplier $\lam \neq 1$; and
\item
when $m \ge 2$, $\al$ is attracting with multiplier $\lam$ 
satisfying $|\lam -1/2|\le 1/2$.
\end{itemize}
In particular, if the simplicity hypothesis also holds, 
all non-trivial zeros of $\zeta$ are indifferent fixed point of $\nu_\zeta$.
\end{prop}
Hence (a) implies (b) in Theorem \ref{thm_zeta}. 
Now we show the converse.

\begin{prop}
\label{prop_indiff_then_RH_SC}
If the fixed points of $\nu_\zeta$ in the critical stripe $\mathcal{S}$
are all indifferent, then both 
the Riemann hypothesis and the simplicity hypothesis are affirmative.
\end{prop}

\parag{Proof.}
Let $\al$ be an indifferent fixed point of 
$\nu_\zeta$ in the critical stripe $\mathcal{S}$.
Since $\zeta$ has no pole in $\mathcal{S}$, 
$\al$ is a zero of some order $m \ge 1$ of $\zeta$ 
and the holomorphic index is 
$\iota(\nu_\zeta, \al)= m\al$
by \propref{prop_nu_index}.
The point $1-\al$ is also a zero of $\zeta$ of order $m$ 
contained in $\mathcal{S}$,
and the holomorphic index is $\iota(\nu_\zeta,1-\al) = m(1-\al)$. 

By assumption, both $\al$ and $1-\al$ are indifferent
 fixed points of $\nu_\zeta$. 
Hence by \propref{prop_classification},
the real parts of  
$\iota(\nu_\zeta, \al)=m\al$
and 
$\iota(\nu_\zeta, 1-\al)=m(1-\al)$
are both $1/2$.
This happens only if $m=1$ and $\Re \al=1/2$. 
\QED

\medskip
Let us finish the proof of \thmref{thm_zeta}:

\parag{Proof of \thmref{thm_zeta}.}
The equivalence of (a) and (b) is shown 
by \propref{prop_RH_implies_indiff_or_att}  
and \propref{prop_indiff_then_RH_SC} 
above. 
The condition (b) implies (c) since the fixed points of $\nu_\zeta$ 
off the critical stripe are all repelling by 
\propref{prop_rep_off_crit_stripe}.

Suppose that (c) holds. 
Then any fixed point $\al$ 
of $\nu_\zeta$ in the critical stripe 
is repelling or indifferent, 
and it is also a zero of $\zeta$ of some order $m$
by \propref{prop_nu_index}.
Hence the holomorphic index $\iota(\nu_\zeta,\al)=m\al$ 
satisfies $\Re m\al \in (0,1/2]$.
Moreover, 
$1-\al$ is also a zero of $\zeta$
with the same order $m$ and
$\iota(\nu_\zeta,1-\al)=m(1-\al)$ 
satisfies $\Re m(1-\al) \in [m-1/2,m)$.

If $m \ge 2$, then $\Re m(1-\al)>1/2$
and thus $1-\al$ is attracting.
This is a contradiction.
If $m=1$ and $\al$ is repelling, 
then $1-\al$ is an attracting fixed point.
This is also a contradiction.
Hence $m=1$ and $\al$ is indifferent.
This implies $\Re \al=1/2$ and 
we conclude that (c) implies (a).

The equivalence of (c) and (d) comes from \propref{prop_att_rep}.
\QED

\parag{Remark.}
We used the functional equation to show the equivalence of (a), (b), and (c).

\begin{figure}[htbp]
\begin{center}
\includegraphics[width=.85\textwidth]
{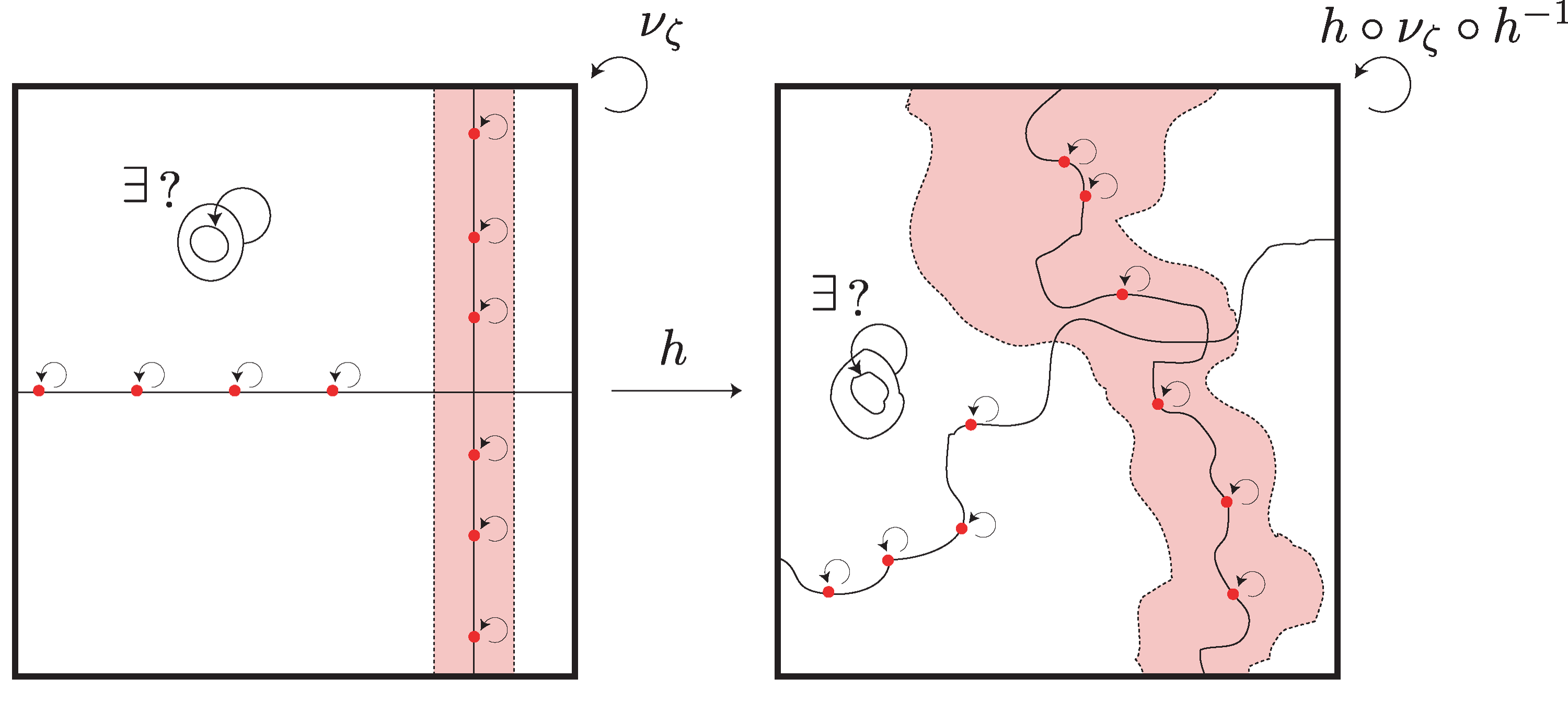}
\end{center}
\caption{Topological deformation of the dynamics of $\nu_\zeta$.
By Theorem \ref{thm_zeta_2},  
we can still state the Riemann hypothesis in the deformed dynamics.}
\label{fig_conj}
\end{figure}

\parag{A more topological version.}
Note again that (d) of \thmref{thm_zeta} is a purely topological 
condition for the function $\nu_\zeta$. 
Even if one observe the dynamics (the action) 
of $\nu_\zeta:\C \to \Chat$ through {\it any} homeomorphism on the sphere,
this condition will be preserved (\figref{fig_conj}). 
More precisely, we have:

\begin{thm}\label{thm_zeta_2}
The conditions {\rm (a) - (d)} of \thmref{thm_zeta} are equivalent to:
 \begin{itemize}
\item[\rm (e)]
For any homeomorphism $h:\Chat \to \Chat$ with $h(\infty)=\infty$, 
the continuous function $\nu_{\zeta,h}:=h \cc \nu_\zeta \cc h^{-1}$ 
has no topological disk $D$ with $\Bar{\nu_{\zeta,h}(D)} \subset D$.
\end{itemize}
\end{thm}
The proof is a routine.
We may regard the dynamics of the map $\nu_{\zeta,h}$ as a topological deformation
of the original dynamics of $\nu_\zeta$. 
Remark that the critical line may not be a ``line" any more 
when it is mapped by a homeomorphism. 
(That may even have a positive area!)

\parag{Variants.}
Next we consider $\nu_g$ for $g = \xi$, 
the {\it Riemann xi function}
given in the first section. 
It is known that $\xi: \C \to \C$ 
is an entire function whose zeros are exactly 
the non-trivial zeros of the Riemann zeta. 
Moreover, we have a functional equation
$ \xi(z) = \xi(1-z)$. 

Now we have a variant of \thmref{thm_zeta}:
\begin{thm}[Interpretation by xi]\label{thm_xi}
The following conditions are equivalent:
\begin{enumerate}[\rm (a')]
\item
The Riemann hypothesis is affirmative and all zeros of $\zeta$ are simple.
\item 
All the fixed points of the meromorphic function 
$\nu_\xi(z) := z-\dfrac{\xi(z)}{z\,\xi'(z)}$
are indifferent.
\item 
The meromorphic function $\nu_\xi$ above has no attracting fixed point.
\item 
There is no topological disk $D$ with $\Bar{\nu_\xi(D)} \subset D$.
\end{enumerate}
\end{thm}
The proof is just analogous to the case of $\nu_\zeta$.
(It is simpler because $\xi$ has neither trivial zeros nor poles.) 

Here is another example. 
To have an entire function it is enough to consider the function 
$$
\eta(z):=(z-1)\zeta(z),
$$
which seems simpler than $\xi(z)$.
In this case $\eta(z)$ and $\zeta(z)$ share all the zeros and 
we will have a similar theorem to \thmref{thm_zeta}.

Yet another interesting entire function is
$$
\chi(z):=(z-a)^m\eta(z)=(z-a)^m(z-1)\zeta(z)
$$ 
where $m \in \N$ and $a \in \C$ satisfy $\Re ma > 1/2$. 
It has an extra zero at $z=a$ of order $m$
 which is an attracting fixed point of $\nu_{\chi}$
by \propref{prop_nu_index}.
The virtue of this family with parameters $m$ and $a$ 
is that we can characterize the chaotic locus (the {\it Julia set} defined in the next section) of the dynamics of $\nu_{\chi}$ as the boundary of the {\it attracting basin} of $a$ (i.e., 
the set of points $z$ whose orbit $z, \nu_\chi(z), 
\nu_\chi(\nu_\chi(z)), \ldots$ converges to the attracting fixed point $a$.) In general, this property helps us to draw the pictures of the chaotic locus. 

\section{Global dynamics and linearization problem}

If any of conditions in \thmref{thm_zeta} is true,
then all the fixed points of $\nu_\zeta$ (and $\nu_\xi$) 
are indifferent. 
In this section we give some brief observations on these fixed points. 

\parag{Global dynamics of $\nu_g$.}
Let $g:\C \to \Chat$ be a non-constant meromorphic function 
and consider the dynamics given by iteration of $g$:
$$
\C \mapsTo{g} \Chat  \mapsTo{g} \Chat  \mapsTo{g} \Chat  \mapsTo{g} 
 \cdots
$$ 
(See \cite{Mi} or \cite{B} for basic results in complex dynamics.)
We use a conventional notation $g^n~(n \ge 0)$ to denote $n$ times
iteration of $g$. 
More precisely, we set $g^0:=\id$ and $g^n:=g \cc g^{n-1}$ for $n \in \N$.
The {\it Fatou set} $F(g)$ of $g$ is defined by the set of $z \in \C$
which has a neighborhood $U$ where $g^n$ is defined for all $n \in \N$ and 
the family 
$\skakko{z \mapsto g^n(z)}_{n \in \N}$ 
is normal (or equivalently, equicontinuous with respect to the 
spherical metric of $\Chat$).
The complement $\Chat -F(g)$ is called 
the {\it Julia set} of $g$ and denoted by $J(g)$. 
We regard the Fatou set and the Julia set as 
the stable and chaotic parts of the dynamics.
For example, the attracting (resp. repelling) 
fixed points (periodic points, more generally)
belong to the Fatou (resp. Julia) set. 

Indifferent fixed points may be in the Fatou set or 
the Julia set, and this is a source of the 
{\it linearization problem}.

\parag{Classification of indifferent fixed point.}
(See \cite[\S 10, \S 11]{Mi} for more details.)
Let $g:\Omega \to \C$ be a holomorphic function
on a domain $\Omega$ and $\al$ 
an indifferent fixed point of multiplier $g'(\al)=e^{2\pi i \theta}~
(\theta \in \R)$. 
In the theory of complex dynamics, 
we say $\al$ is {\it parabolic}
or {\it rationally indifferent}
 if $\theta \in \Q$,
and {\it irrationally indifferent} otherwise. 

Let $\al$ be an irrationally indifferent fixed point of $g$
with $g'(\al)=e^{2\pi i \theta}~(\theta \in \R)$.
It is called {\it linearizable} if there exists a holomorphic 
homeomorphism $\phi$ defined on a neighborhood $U$ of $\al$ 
such that $\phi(\al)=0$; $\phi(U)$ is a round disk centered at $0$; and 
$\phi \cc g \cc \phi^{-1}(w)=e^{2 \pi i \theta} w$.
In other words, the local action near $a$ is
 conjugate to a rotation of angle $2 \pi \theta$.
If this is the case, the largest possible neighborhood 
$U$ with these properties is called a {\it Siegel disk}.
The dynamics of $g$ on the Siegel disk 
is equicontinuous and thus $U$ is contained in the Fatou set.
In fact, it is known that an indifferent fixed point 
belongs to the Fatou set if
and only if it is linearizable.

It is also known that for Lebesgue almost every $\theta \in [0,1) -\Q$, 
the irrationally indifferent fixed point of multiplier $e^{2\pi i \theta}$ is
linearizable. 
For example, suppose that $\theta \in [0,1)-\Q$ 
has the continued fraction expansion 
$$
\theta = \frac{1}{a_1+\dfrac{1}{a_2+\dfrac{1}{a_3+\cdots}}}
=:[a_1,\,a_2,\,a_3,\,\ldots],
$$
where the coefficients $a_i \in \N$ are uniformly bounded.
Such a $\theta$ is called of {\it bounded type}. 
Then the irrationally indifferent fixed points of 
multiplier $e^{2\pi i\theta}$ with bounded type 
$\theta$ is always linearizable.
However, there is also a dense subset of 
$\theta \in [0,1) -\Q$
that gives non-linearizable fixed points.
(Unfortunately, the multiplier $e^{2 \pi i \theta}$ 
is not the only factor for the precise condition of linearizability.)

\parag{Does $\nu_\zeta$ has Siegel disks?}
Under the Riemann hypothesis and the simplicity hypothesis,
each non-trivial zero of $\zeta$ is 
of the form $\al=1/2+\gam i~(\gam \in \R)$ 
and $\al$ is an indifferent fixed point 
of $\nu_\zeta$ (or $\nu_\xi$, etc.).

Now it is easy to check:
\begin{prop}
\label{prop_rotation}
Under the Riemann hypothesis and the simplicity hypothesis,
non-trivial zero $\al=1/2+\gam i~(\gam \in \R)$
is an indifferent fixed point of $\nu_\zeta$ 
with multiplier $e^{2\pi i \theta}$
where the values $\gam$ and $\theta$ are related by
$$
\gam \ee \frac{1}{2 \tan \pi \theta} 
~~\Longleftrightarrow~~ \theta \ee \frac{1}{\pi}\arctan \frac{1}{2\gam}.
$$
\end{prop}
Note that $\theta \to 0$ as $\gam \to \pm \infty$.

\parag{Proof.}
By \propref{prop_nu_index},
we have $\nu_\zeta'(\al) = 1-1/\al$ and thus
$e^{2\pi i \theta} = 1-1/(1/2 + \gam i)$.
This is equivalent to the formulae above.  
\QED

\medskip
Now it is natural to ask:
\begin{quote}
\parag{Linearization problem.}
{\it 
Can $\theta$ be a rational number?
Is $\al$ linearizable? 
That is, can $\nu_\zeta$ has a Siegel disk?
}
\end{quote}

\parag{Numerical observation.}
As we have mentioned, randomly chosen 
$\theta$ gives a Siegel disk with probability one.
Table \ref{table_rotation_numbers} gives 
the continued fractions up to the 50th term for some zeros. 
(Here $1/2 + \gam_n i~(\gam_n>0)$ is the $n$th non-trivial 
zero of $\zeta$ from below and 
$\gam_n = 1/(2\tan \pi \theta_n)$.)
It seems very unlikely that these zeros have
rational $\theta$.

{\small
\begin{table}[htbp]
\begin{center}
\begin{tabular}{|l|l|l|l|}\hline
$n$ & $\Im \gam_n$ & $\theta_n$ & $[a_1,a_2,\cdots,a_{25},a_{26},\cdots, a_{50}]$
\\\hline
1 &14.1347 & 0.0112552 & 
[88,1,5,1,1,2,2,5,2,15,2,4,2,4,1,9,1,1,5,2,10,1,1,5,1,
\\ & & &
2,7,100,9,2,2,3,2,5,1,8,179,23,1,1,35,1,3,1,2,8,7,34,4,1]
\\\hline
2 &21.0220 & 0.00756943 & 
[132,9,14,1,1,1,2,1,52,1,9,3,4,1,1,1,1,2,2,3,2,1,10,1,1,
\\ & & &
1,9,1,1,6,5,1,5,3,1,5,2,6,1,135,1,1,5,1,2,3,2,9,1,3]
\\\hline
3 &25.0109 & 0.00636259 &
[157,5,1,12,3,1,1,1,1,2,11,1,29,5,1,4,1,1,3,5,14,1,3,1,2,
\\ & & &
3,6,14,4,1,41,4,1,1,7,4,1,3,21,1,4,3,1,2,2,1,16,1,2,3]
\\\hline
4 & 30.4249 & 0.00523061 & 
[191,5,2,15,3,2,2,7,2,1,2,46,2,1,1,6,1,4,2,2,4,1,6,1,1,
\\ & & &
2,5,1,8,1,2,2,5,1,4,39,3,19,5,2,9,1,1876,2,12,1,4,4,1,6]
\\\hline\hline
10 & 49.7738 & 0.00319746 &
[312,1,2,1,48,1,1,4,1,3,1,2,5,1,21,1,4,1,3,2,1,1,8,1,6,
\\ & & &
9,1,2,3,2,1,2,1,1,1,4,1,1,1,5,2,1,3,126,1,24,3,2,29,5]
\\\hline
$10^2$ & 236.524 & 0.00067289 &
[1486, 7,1,4,1,2,1,53,2,8,1,4,6,3,1,3,1,13,3,1,7,2,18,3,1,
\\ & & &
1,1,1,1,2,1,210,3,4,1,3,1,2,2,8,1,1,7,2,2,2,1,12,2,2]
\\\hline
$10^3$ & 1419.42 & 0.000112127 &
[8918,2,48,220,1,15,1,1,6,1,1,6,2,4,149,1,15,3,6,1,4,1,1,64,11,
\\ & & &
1,1,13,11,2,4,2,3,2,4,1,2,1,1,2,3,46,1,11,3,1,18,1,6,4]
\\\hline
$10^4$ & 9877.78 & 0.0000161124 & 
[62063,1,15,2,1,2,1,8,1,2,6,2,1,2,2,1,5,1,1,186,3,5,4,14,9,
\\ & & &
1,12,1,12,1,6,4,22,1,3,15,1,2,1,14,1,2,2,12,3,1,1,1,45,2]
\\\hline
\end{tabular}

\end{center}
\caption{Continued fraction expansions for some zeros of $\zeta$.}
\label{table_rotation_numbers}
\end{table}
}

\begin{figure}[htbp]
\begin{center}
\includegraphics[height = .85\textheight, bb = 0 0 768 2765]{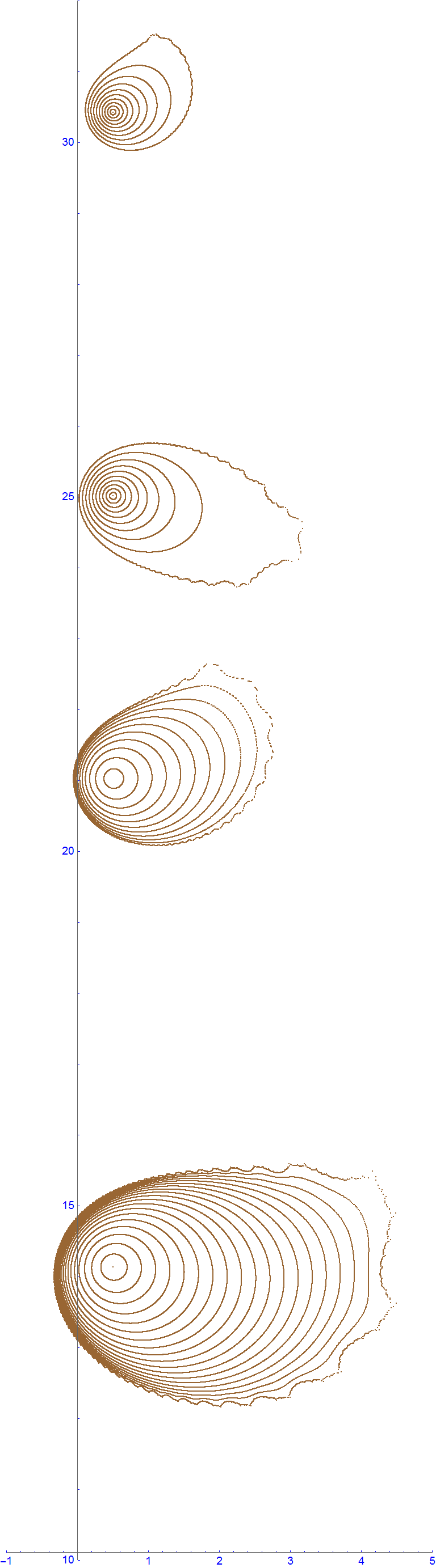}~~~~
\includegraphics[height = .85\textheight, bb = 0 0 256 2458]{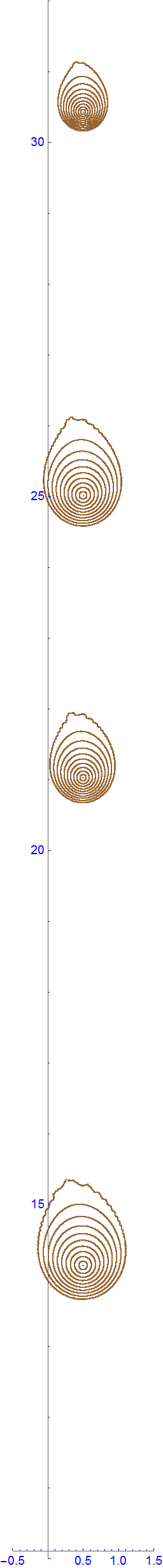}
\end{center}
\caption{The first four Siegel disks (?) of $\nu_\zeta$ (left) and $\nu_\xi$ (right).}
\label{fig_siegels}
\end{figure}

\figref{fig_siegels} shows some recurrent orbits near 
the first 4 zeros of $\zeta$ in the dynamics of $\nu_\zeta$ and $\nu_\xi$.

\section{Appendix: Newton's method}
There are many root finding algorithms, 
but the most famous one would be \textit{Newton's method}. 
Let us apply it to the Riemann zeta and its variants.
The aim of this Appendix is to describe it 
in terms of holomorphic index.

\parag{Relaxed Newton maps and fixed points.}
Let $g:\C \to \Chat$ be a non-constant 
meromorphic function and $\kappa$ a complex constant
with $|\kappa-1|<1$. 
We define its {\it relaxed Newton's map} $N_g=N_{g,\kappa}: \C \to \Chat$ by
$$
N_g(z) \ee z - \kappa \, \frac{g(z)}{g'(z)},
$$
which is again meromorphic. (See \cite{B}.) 
When $\kappa=1$ the map $N_g=N_{g,1}$ 
is traditional {\it Newton's map} for Newton's method.
Here is a version of \propref{prop_nu_index}
(and \corref{cor_no_fixed_points_with_mult_one})
for $N_g$:

\begin{prop}[Fixed points of $N_g$]
\label{prop_newton_index}
The point $\al \in \C$ is a fixed point of $N_g$
if and only if 
$\al$ is a zero or a pole of $g$. 
Moreover, 
\begin{itemize}
\item
if $\al$ is a zero of $g$ of order $m \ge 1$,
then $\al$ is an attracting fixed point of
multiplier $N_g'(\al) =1 - \kappa/m$  
and its index is $\iota(N_g,\al)= m/\kappa$. 
\item 
if $\al$ is a pole of $g$ of order $m \ge 1$,
then $\al$ is a repelling fixed point of
multiplier $N_g'(\al) =1 +\kappa/m$  
and its index is $\iota(N_g,\al)= -m/\kappa$.
\end{itemize}
In particular, $N_g$ has no fixed point of multiplier $1$.
\end{prop}

The proof is similar to 
that of \propref{prop_nu_index} 
and left to the readers.

\parag{Newton's method.}
The idea of (relaxed) Newton's method is 
to use the attracting fixed points
 of $N_g$ to detect the location of the zeros of $g$. 
More precisely, by taking an initial value $z_0$ 
sufficiently close to the zero $\al$,
the sequence $\skakko{N_g^n(z_0)}_{n \ge 0}$ 
converges rapidly to the attracting fixed point $\al$. 
It is practical to use the traditional value 
$\kappa=1$, since the convergence 
to the simple zero is quadratic. 
That is, we have 
$N_g(z)-\al=O((z-\al)^2)$ near $\al$.

\parag{Holomorphic index and the argument principle.}
Let us restrict our attention to  the traditional case $\kappa=1$.
For the fixed point $\al$ of $N_g$ its holomorphic index is
$$
\iota(N_g,\al)
=\frac{1}{2\pi i}\int_C \frac{1}{z-N_g(z)}\dz
=\frac{1}{2\pi i}\int_C \frac{g'(z)}{g(z)}\dz
$$
where $C$ is a small circle around $\al$.
This is exactly the argument principle 
applied to $g$. 

\parag{The Riemann hypothesis.}
When we apply Newton's method to the Riemann zeta,
all zeros of $\zeta$ become attracting fixed points of 
$N_\zeta(z) = z-\zeta(z)/\zeta'(z)$.
By \propref{prop_att_rep}, we have:

\begin{thm}
\label{thm_NZeta}
The Riemann hypothesis is affirmative
if and only if there is no topological disk $D$ 
contained in the stripe $\mathcal{S}'  
= \skakko{z \in \C \st 1/2<\Re z <1}$
satisfying $\overline{N_\zeta(D)} \subset D$. 
\end{thm}
 
\parag{Proof.}
The Riemann hypothesis holds if and only if 
there is no zero of $\zeta$ in the 
stripe $\mathcal{S}' \subset \mathcal{S}$.
Equivalently, by \propref{prop_newton_index}, 
there is no attracting fixed point of 
$N_\zeta$ in $\mathcal{S}'$. 
Now the theorem is a simple application of \propref{prop_att_rep}.
\QED

\parag{Some pictures.} 
It is easier to draw pictures of the Julia set 
of Newton's map $N_\zeta$ than those of $\nu_\zeta$.
After the list of references, we present some pictures of $N_\zeta$, $N_\eta$, $N_\xi$ and $N_{\cosh}$ 
with comments.

\begin{figure}[h]
\begin{center}
\includegraphics[width = .475\textwidth,bb = 0 0 800 1067]{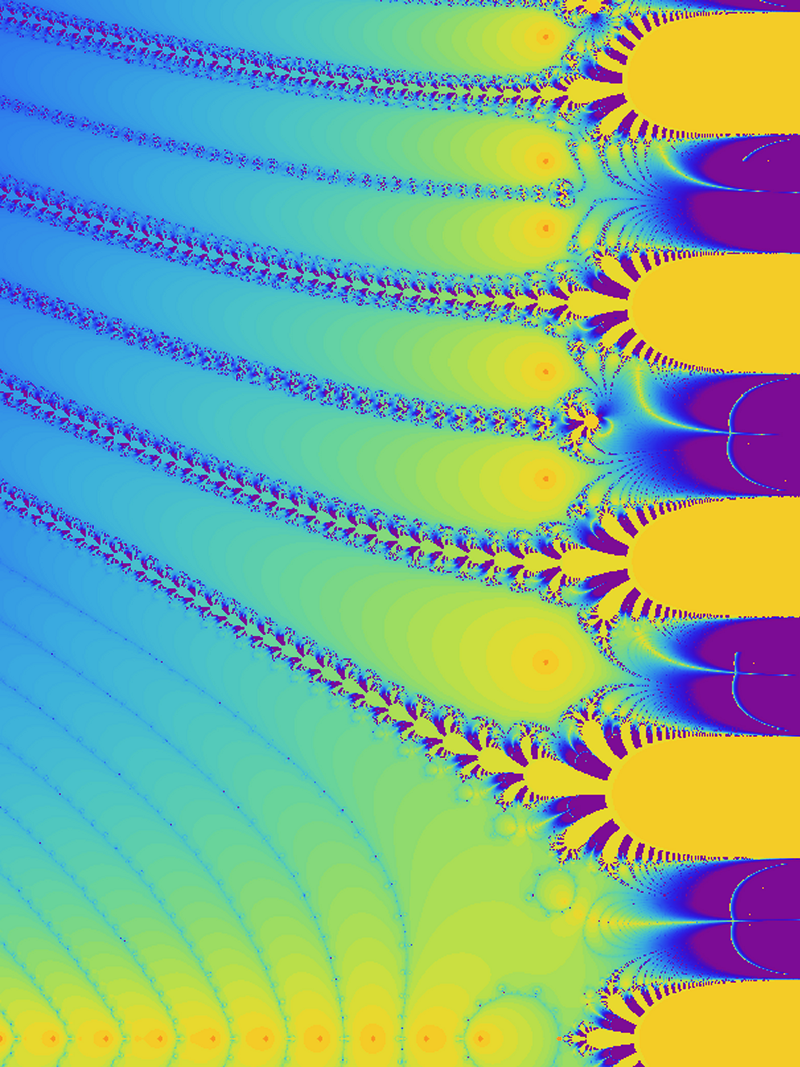}
\includegraphics[width = .475\textwidth,bb = 0 0 800 1067]{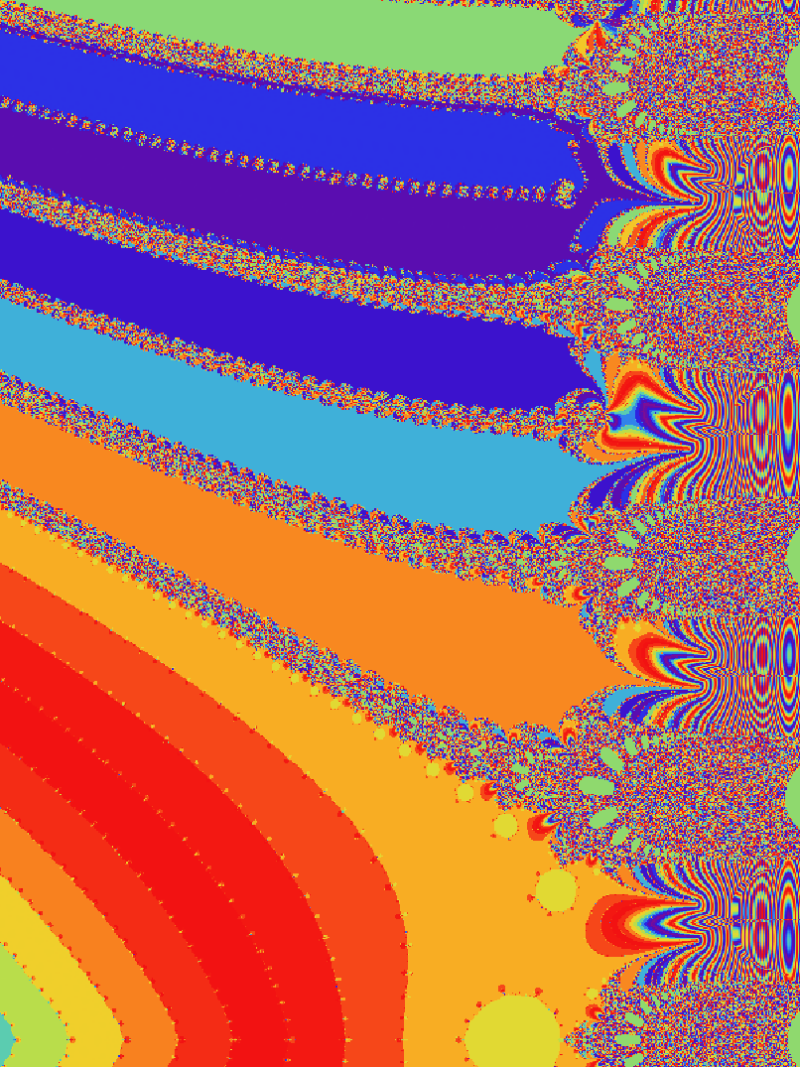}
\end{center}
\caption{Left: The Julia set of $N_\zeta$ 
in $\{\Re z \in [-20,10],\, \Im z \in [-1,39]\}$. The orange dots 
are disks that are close to the zeros of $\zeta$. 
Right: The same region in different colors. 
Colors distinguish the zeros to converge.
}
\label{fig_NZeta}
\end{figure}

\begin{figure}[h]
\begin{center}
\includegraphics[width = .6\textwidth,bb = 0 0 806 807]{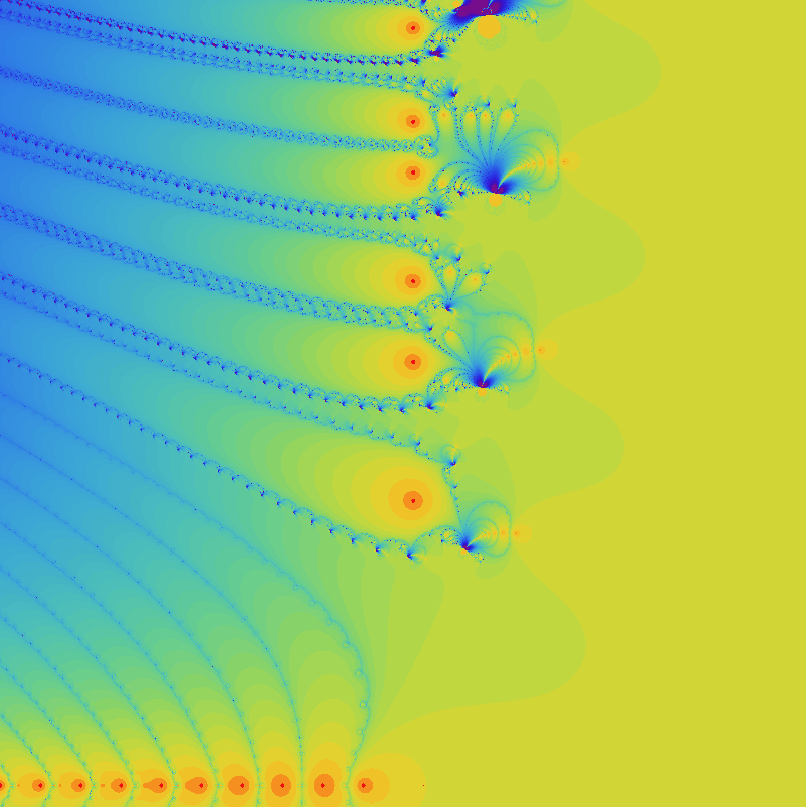}\\[.5em]
\includegraphics[width = .6\textwidth,bb = 0 0 806 807]{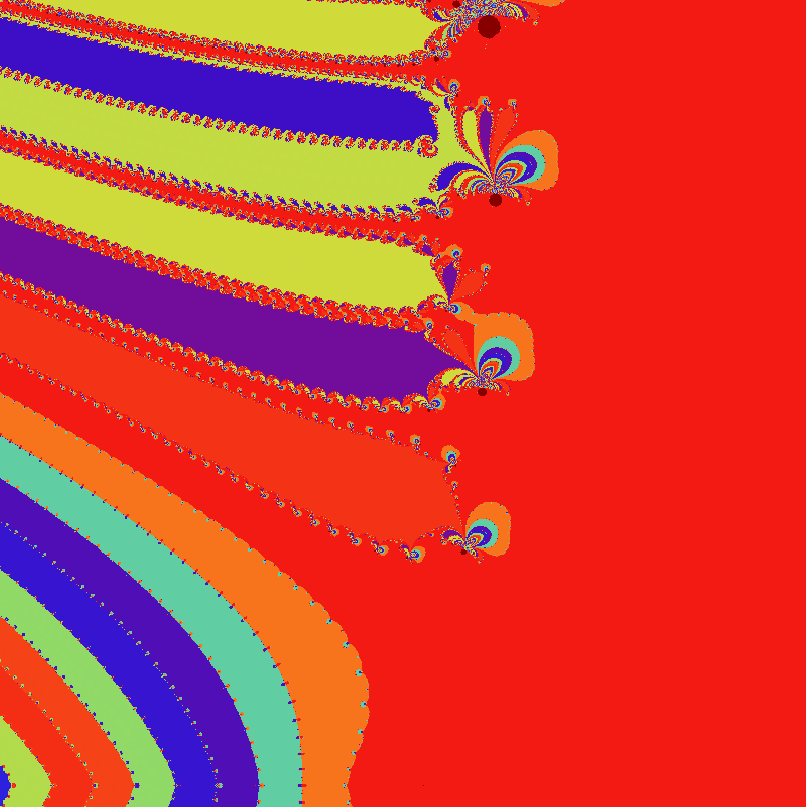}
\end{center}
\caption{
The Julia set of $N_\eta$ in $\{\Re z \in [-20,20],\, \Im z \in [-1,39]\}$,
drawn in the same colors as Figure \ref{fig_NZeta}. 
Probably because $\eta$ is entire, 
the dynamics of $N_\eta$ is simpler than that of $N_\zeta$. 
It is known that for any entire function $g$ and its zeros, 
their immediate basins (the connected components of the Fatou set of $N_g$ 
which contain the zeros) are simply connected and unbounded (\cite{MS}). 
}
\label{fig_NEta}
\end{figure}

\begin{figure}[h]
\begin{center}
\includegraphics[width = .3\textwidth,bb = 0 0 278 1160]{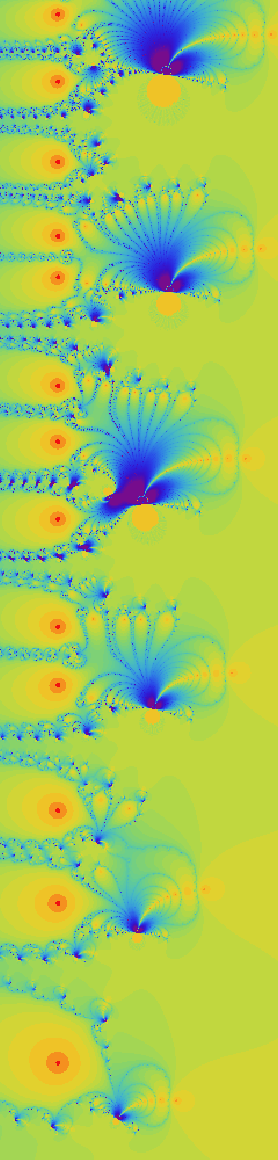}
\includegraphics[width = .3\textwidth,bb = 0 0 278 1160]{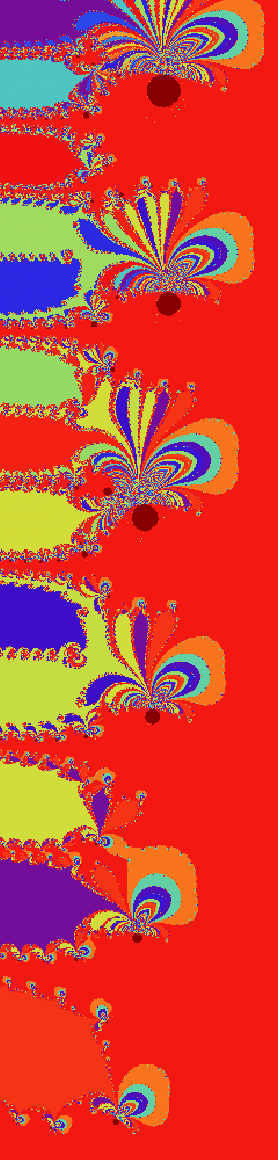}
\end{center}
\caption{
Details of the Julia set of $N_\eta$ 
in $\{\Re z \in [-2,10],\, \Im z \in [10,60]\}$ (``Heads of Chickens").
}
\label{fig_NEta-chickens}
\end{figure}

\begin{figure}[h]
\begin{center}
\includegraphics[width = .475\textwidth,bb = 0 0 800 800]{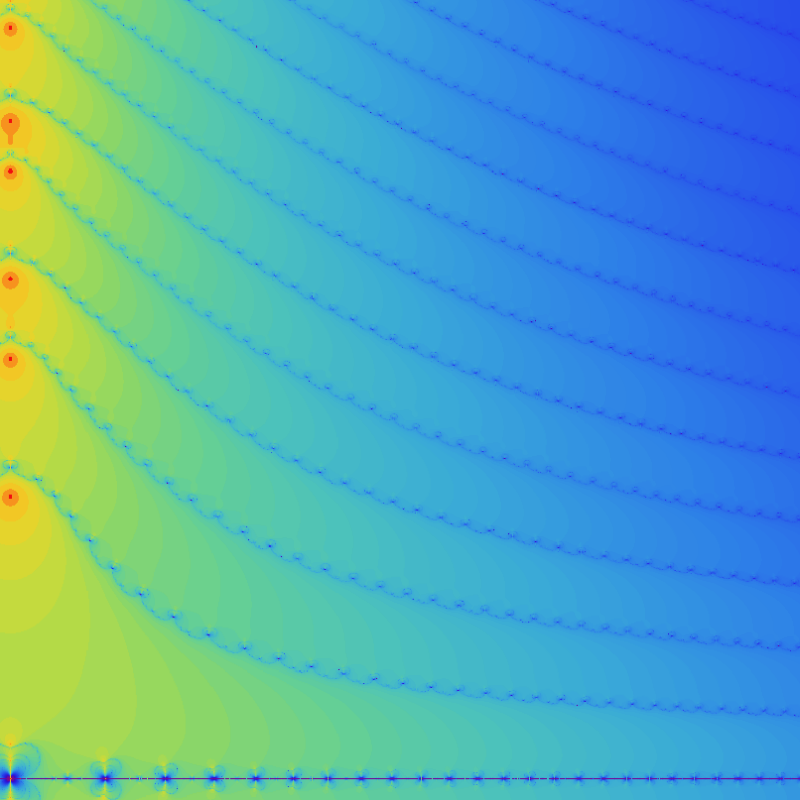}
\includegraphics[width = .475\textwidth,bb = 0 0 800 800]{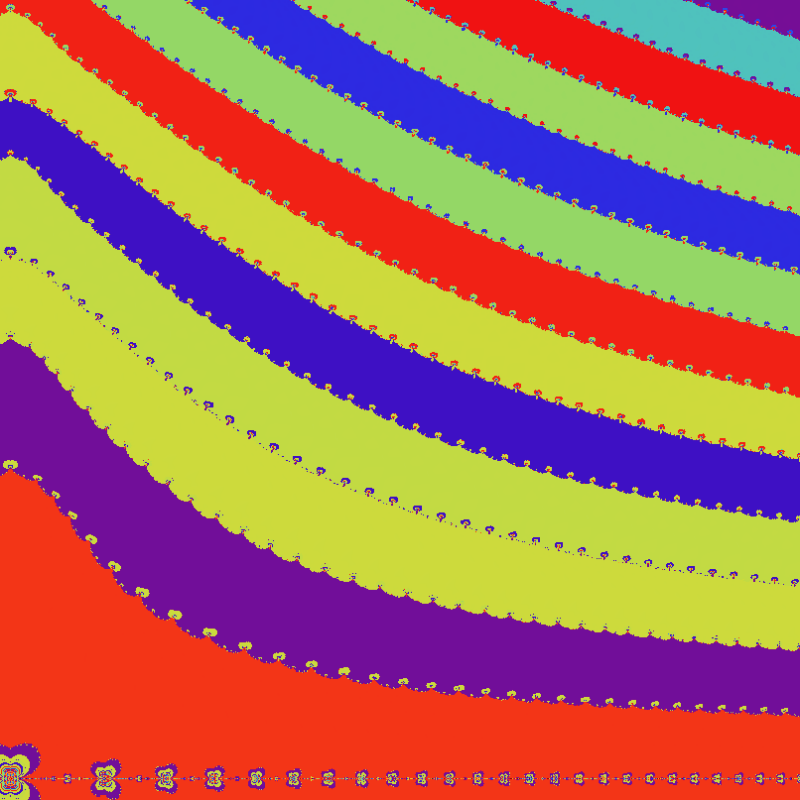}
\end{center}
\caption{
Details of the Julia set of $N_\xi$ 
in $\{\Re z \in [0,40],\, \Im z \in [-1,39]\}$.
The dynamics seems surprisingly simple. 
Compare with the case of hyperbolic cosine in Figure \ref{fig_NCosh}. 
}
\label{fig_NXi}
\end{figure}

\begin{figure}[h]
\begin{center}
\includegraphics[width = .8\textwidth,bb = 0 0 806 336]{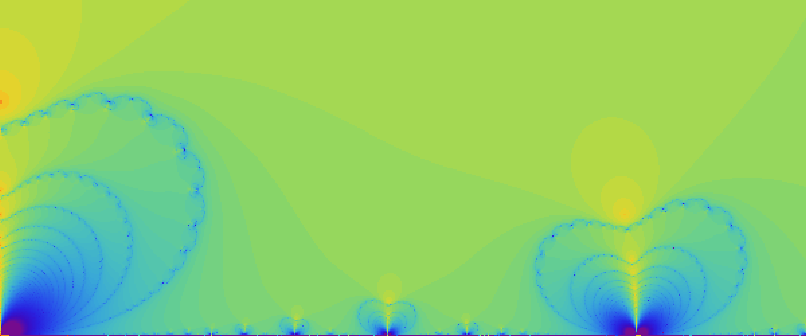}\\
\includegraphics[width = .8\textwidth,bb = 0 0 806 336]{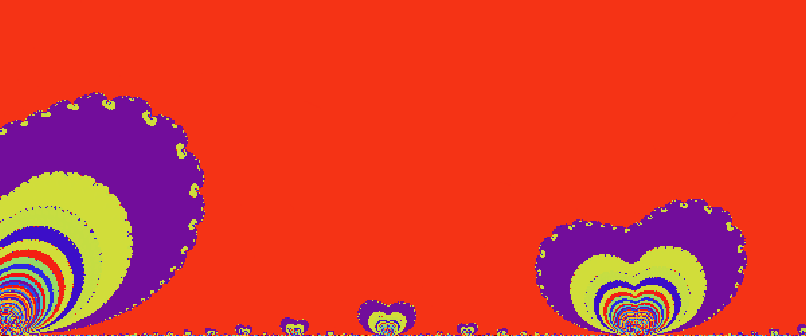}
\end{center}
\caption{Details of the Julia set of $N_\xi$ 
in $\{\Re z \in [0.5,6.5],\, \Im z \in [0,2.5]\}$.
}
\label{fig_NXi-details}
\end{figure}

\begin{figure}[h]
\begin{center}
\includegraphics[width = .475\textwidth,bb = 0 0 806 807]{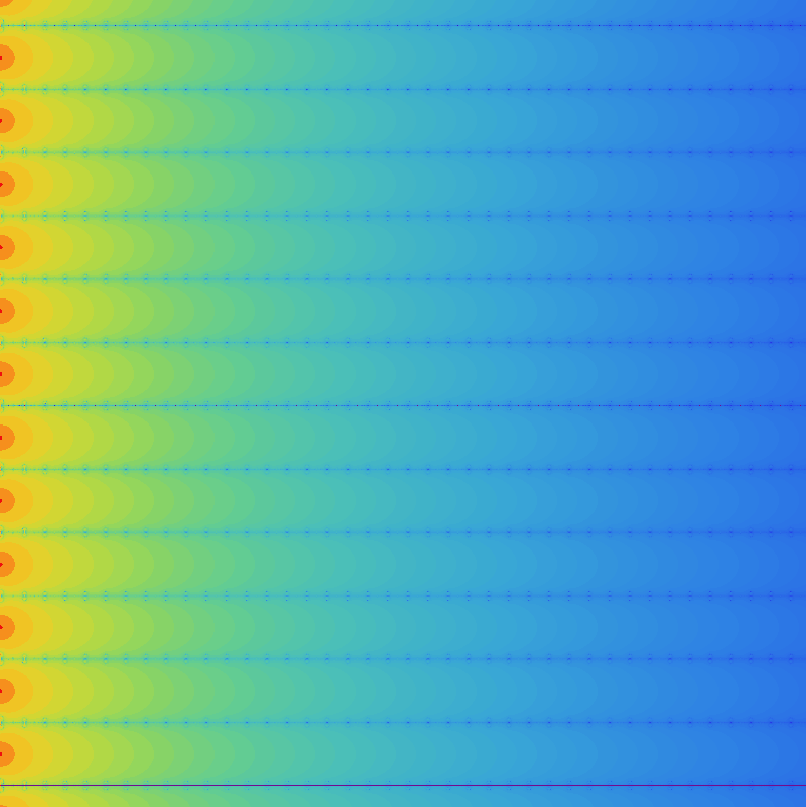}
\includegraphics[width = .475\textwidth,bb = 0 0 806 807]{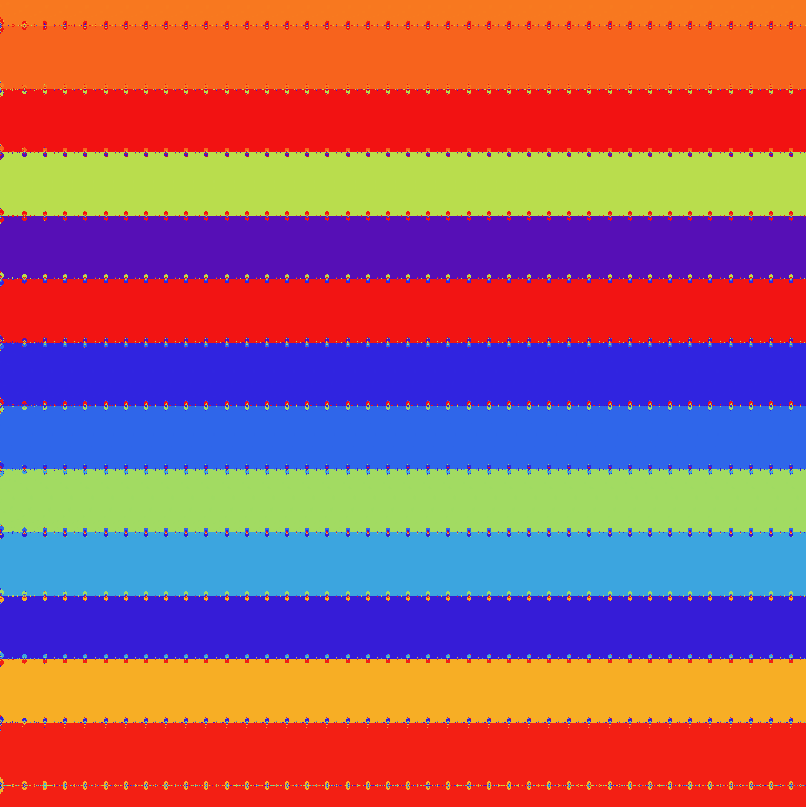}
\end{center}
\caption{
The Julia set of the Newton map $N_{\cosh}$ 
of the hyperbolic cosine $\cosh z$ 
in $\{\Re z \in [0,40],\, \Im z \in [-1,39]\}$.
The Julia set seems topologically the same as that of $N_\xi$.
Do they belong to the same deformation space?
}
\label{fig_NCosh}
\end{figure}

\begin{figure}[h]
\begin{center}
\includegraphics[width = .25\textwidth,bb = 0 0 806 1882]{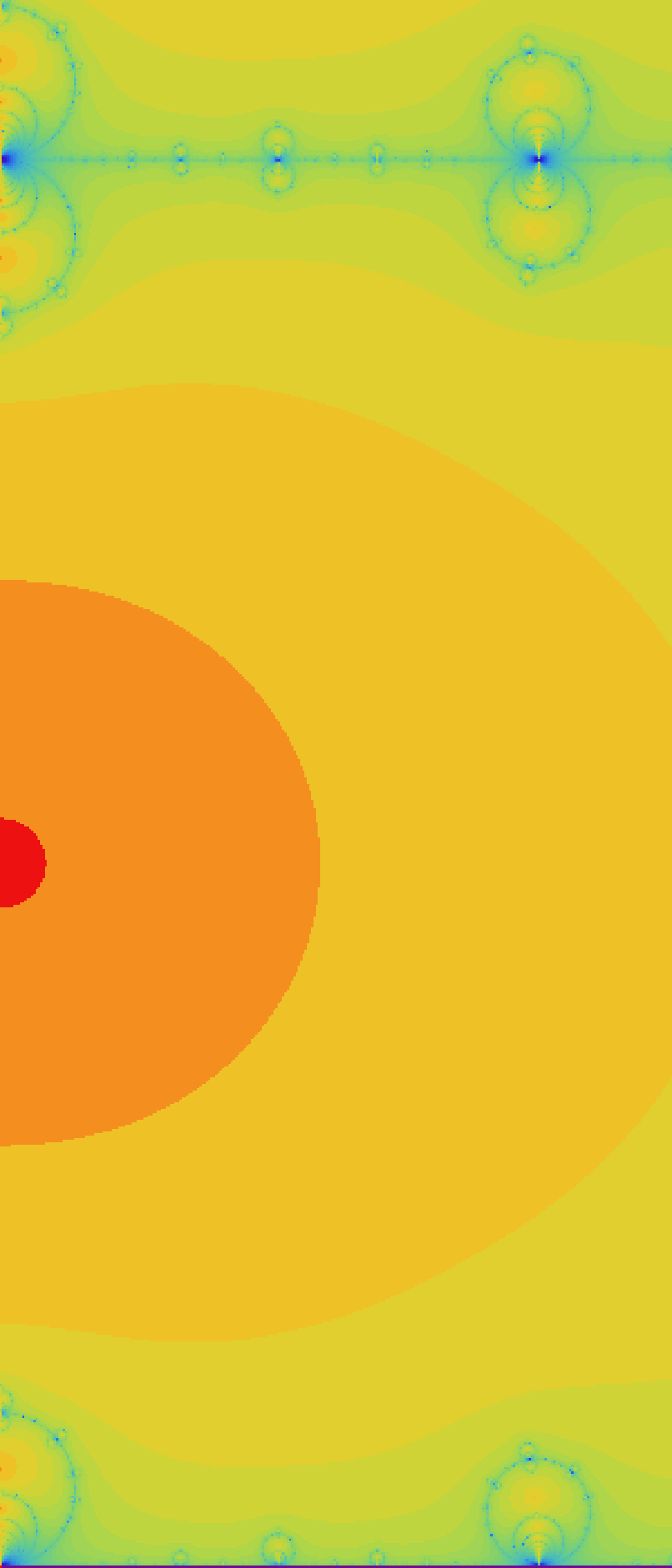}
\includegraphics[width = .25\textwidth,bb = 0 0 806 1882]{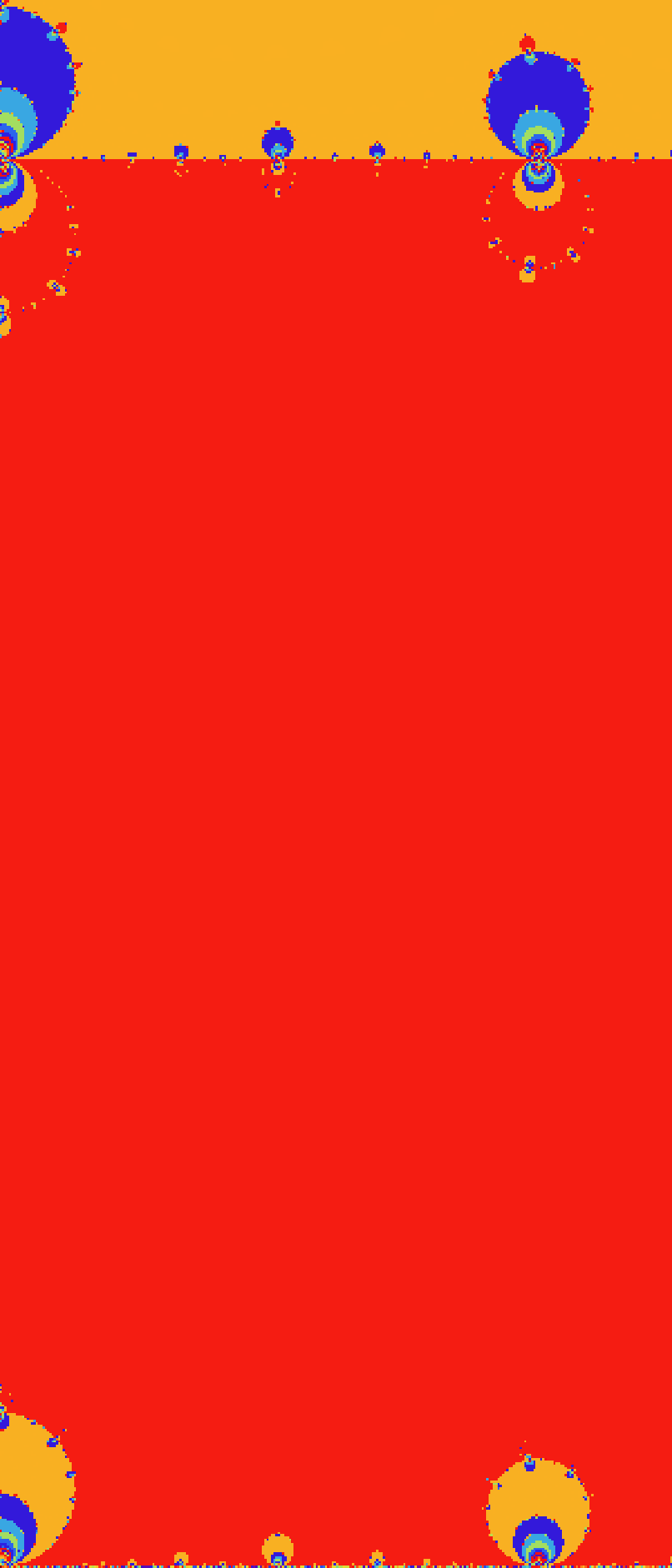}
\end{center}
\caption{
Details of the Julia set of the Newton map $N_{\cosh}$ 
in $\{\Re z \in [0,1.5],\, \Im z \in [0,3.5]\}$.
}
\label{fig_NCosh-details}
\end{figure}

\end{document}